\documentclass[12pt,a4paper]{article}
\usepackage[T2A]{fontenc}
\usepackage[cp1251]{inputenc}
\usepackage[english,russian]{babel}

\usepackage{amsmath}
\usepackage{amsfonts}
\usepackage{amssymb}

\title{Technical report: A generating funtion for the Euler numbers of the second kind and it's application}
\author{Dmitry V. Kruchinin, Vladimir V. Kruchinin}

\begin{document}

\newcommand{\stirlingone}[2]{\left[{#1 \atop #2}\right]}
\newcommand{\stirlingtwo}[2]{\left\{{#1 \atop #2}\right\}}
\newcommand{\eulereantwo}[2]{\left < \!\! \left < {#1 \atop #2} \right > \!\!\right >}

\maketitle

\begin{abstract}
In the paper, 2 explicit formulas for the Euler numbers of the second kind are obtained. Based on those formulas a exponential generating function is deduced. Using the generating function some well-known and new identities for the Euler number of the second kind are obtained.
\end{abstract}

\section{Introduction}

The Euler numbers of the second kind is very useful in combinatorics, number theory, graph theory, analytic geometry and other areas \cite{Gessel,Knuth,Forest}. The Euler numbers of the second kind is a numerical triangle. In the paper, we consider the triangle with an initianal term $\eulereantwo{1}{1}$. This variant of the triangle is submitted in the On-Line Encyclopedia of Integer Sequences by number of secuence A008517\cite{OEIS}.

The triange of the Euler numbers of the second kind is defined by recurrence expression
\begin{equation}\label{To_3}
\eulereantwo{n}{m}= m\,\eulereantwo{n-1}{m}+(2\,n-m)\eulereantwo{n-1}{m-1}.
\end{equation}

A detailed review of the Euler numbers of the second kind was given by Knuth \cite{Knuth}. However, in the present time there is not a generating funtion for those numbers and there are a few explicit ways to defined those numbers.  In this paper we give 2 explicit formulas and a generating function for the Euler numbers of the second kind. Based of those results we obtain new identities for the Euler numbers of the second kind.

\section{Explicit formulas}

{\bf Theorem} For the Euler numbers of the second kind there hold the following formulas:
\begin{equation}\label{Fo_E1}
\eulereantwo{n}{m}=\sum_{k=0}^m (-1)^{m-k}\,{2\,n+1 \choose m-k} \stirlingtwo{n+k}{k},
\end{equation}

\begin{equation}\label{Fo_E2}
\eulereantwo{n}{m}=\sum_{k=0}^m (-1)^{m-k}\,k\,{2\,n \choose m-k}\stirlingtwo{n+k-1}{k}.
\end{equation}

{\bf Proof} First we state the well-known results that we will use:\\
1) For the binomial coefficients we have\\
\begin{equation}\label{To_1}
{n \choose m} = {n-1 \choose m} + {n-1 \choose m-1},
\end{equation}
2) For the Stirling numbers of the second kind we have\\
\begin{equation}\label{To_2}
\stirlingtwo{n}{m}= m\,\stirlingtwo{n-1}{m}+\stirlingtwo{n-1}{m-1}.
\end{equation}

For the proving the formula  (\ref{Fo_E1}) we apply the mathmatical induction. For  $n=1$ and $m=1$ both of formulas are equal to 1. Suppose the formula (\ref{Fo_E1}) is right for every $n$ and $m-1$.
Using the identity (\ref{To_3})  for (\ref{Fo_E1}), we get
\begin{eqnarray}\nonumber
\eulereantwo{n}{m}=m\sum_{k=0}^m (-1)^{m-k}\,{2\,n-1 \choose m-k} \stirlingtwo{n+k-1}{k}+\\
\nonumber
+(2\,n-m)\,\sum_{k=0}^{m-1} (-1)^{m-k-1}\,{2\,n-1 \choose m-k-1} \stirlingtwo{n+k-1}{k}.
\end{eqnarray}

We note that the second kind of the expression for $k=m$ is equal to 0. Then we combine both of sums:
\begin{eqnarray}\nonumber
\eulereantwo{n}{m}=\sum_{k=0}^m (-1)^{m-k}\left(m{2\,n-1 \choose m-k} -(2\,n-m){2\,n-1 \choose m-k-1} \right)\stirlingtwo{n+k-1}{k}.
\end{eqnarray}
Next we consider the difference between coefficients inside the brackets.
$$
m{2\,n-1 \choose m-k} -(2\,n-m){2\,n-1 \choose m-k-1} =
$$
$$
=m{2\,n-1 \choose m-k}+m{2\,n-1 \choose m-k-1}-2\,n{2\,n-1 \choose m-k-1} =
$$
$$
=m{2\,n \choose m-k}-(m-k){2\,n \choose m-k} =k{2\,n \choose m-k}.
$$
Hence, we get the formula (\ref{Fo_E2})
$$
\eulereantwo{n}{m}=\sum_{k=0}^m (-1)^{m-k}\,k\,{2\,n \choose m-k}\stirlingtwo{n+k-1}{k}.
$$
Applying the identity for the Stirling number of the second kind for the formula (\ref{Fo_E1}) we get 
\begin{eqnarray}\label{sum1}
\eulereantwo{n}{m}=\sum_{k=0}^m (-1)^{m-k}{2\,n+1 \choose m-k}\,k\,\stirlingtwo{n+k-1}{k}+\\
+\sum_{k=0}^{m-1} (-1)^{m-k}{2\,n+1 \choose m-k}\stirlingtwo{n+k-1}{k-1}
\end{eqnarray}

Consider the second sum in the above expression we note that for $k=0$ the expression inside sum is equal to 0, because of the Stirling numbers of the second kind have negative second parameter.
$$
\sum_{k=1}^m (-1)^{m-k}{2\,n+1 \choose m-k}\stirlingtwo{n+k-1}{k-1}.
$$
Substituting $k+1$ for $k$, we get
$$
\sum_{k=0}^{m-1} (-1)^{m-k-1}{2\,n+1 \choose m-k-1}\stirlingtwo{n+k}{k}
$$
Therefore, we obtain the formula for $\eulereantwo{n}{m-1}$.

Next we consider sum in (\ref{sum1}) and apply the identity for the binomial coefficients
$$
\sum_{k=0}^m (-1)^{m-k}{2\,n+1 \choose m-k}\,k\,\stirlingtwo{n+k-1}{k}=
$$
$$
=\sum_{k=0}^m (-1)^{m-k}{2\,n \choose m-k}\,k\,\stirlingtwo{n+k-1}{k}+
$$
$$
+\sum_{k=0}^m (-1)^{m-k}{2\,n \choose m-k-1}\,k\,\stirlingtwo{n+k-1}{k}.
$$

The first sum is the formula (\ref{Fo_E2}) that obtained based on the identity for the Euler numbers of the second kind. 
For $k=m$ the expression inside the second sum is equal to 0, because the binomial coefficients is equal to 0. Since that we have
$$
\sum_{k=0}^{m-1} (-1)^{m-k}{2\,n \choose m-k-1}\,k\,\stirlingtwo{n+k-1}{k}
$$
multiplying on $(-1)$ we get
$$
-\sum_{k=0}^{m-1} (-1)^{m-k-1}{2\,n \choose m-k-1}\,k\,\stirlingtwo{n+k-1}{k}.
$$ 
Using (\ref{Fo_E2}), we obtain
$$
\sum_{k=0}^{m-1} (-1)^{m-k-1}{2\,n \choose m-k-1}\,k\,\stirlingtwo{n+k-1}{k}=-\eulereantwo{n}{m-1}.
$$
Since that, we arrive to desired result
$$
\sum_{k=0}^m (-1)^{m-k}\,{2\,n+1 \choose m-k} \stirlingtwo{n+k}{k}=\eulereantwo{n}{m}-\eulereantwo{n}{m-1}+\eulereantwo{n}{m-1}=\eulereantwo{n}{m}.
$$

\section{Generating function for the Euler numbers of the second kind}

Using the obtained explicit formulas, now we can find an exponential generating function for the Euler numbers of the second kind.\\ 
{\bf Theorem 2} The exponential generating function for the Euler numbers of the second kind is defined by the following expression:
\begin{equation}\label{gfEu}
\sum\limits_{n \geqslant 0} \sum_{m\geqslant \geqslant0} \eulereantwo{n}{m}\frac{x^n}{n!}t^m={{1-t}\over{W\left(-t\,e^{\left(1-t\right)^2\,x-t}
 \right)+1}}, 
\end{equation}
where $W(x)$ -- Lambert function and $x>0$.
	
{\bf Proof}	Present the generating function as dual series with shift $(n-1)$
$$
E(x,t) = \sum_{n>0}\sum_{m>0} \eulereantwo{n-1}{m}\frac{x^n}{n!}t^m 
$$
Preset the formula(\ref{Fo_E1}) as dual series with shift as follows
$$
E(x,t)=\sum_{n>0}P_n(t)\frac{x^n}{n!}\,(1-t)^{2n-1},
$$
where
$$
(1-t)^{2\,n-1}=\sum_{n\geqslant 0}(-1)^m\,{2\,n-1 \choose m}t^m
$$
and
$$
P_n(t)=\sum_{m\geqslant 0}\stirlingtwo{n+m-1}{m}t^m
$$ 
Define the following generating function
$$
u(x,t)=\sum_{n>0}\sum_{m\geqslant 0}\stirlingtwo{n+m-1}{m}\frac{x^n}{n!}t^m
$$
and show that this generating function is a compositional inverse generating function for
$$
y(x,t)=(x-t(e^x-1)).
$$

Suppose the reciprocal generating function for $y(x,t)$ is 
$$
\frac{x}{y(x,t)}=\frac{1}{1-t\frac{e^x-1}{x}}.
$$
Then we find an expression for $k$ powers. For that we write
$$
\left(t\left(\frac{e^x-1}{x}\right)\right)^k=\sum_{n>0} \sum_{m>0}T(n,m,k)x^nt^n
$$
Since for $(e^x-1)$ there is the following identity
$$
(e^x-1)^k=\sum\limits_{n\geqslant 0} k!\stirlingtwo{n}{k}\frac{x^n}{n!},
$$
we get
$$
T(n,m,k)=\delta(m,k) \stirlingtwo{n+k}{k}\frac{k!}{(n+k)!},
$$
where $\delta(m,k)$ is the Kroneker symbol.

Since for  $\frac{1}{(1-x)}$ there is the following identity
$$
\frac{1}{(1-x)^k}=\sum_{n\geqslant 0} {n+k-1 \choose n}x^n.
$$
Then for $k$ powers of the composition of the two generating function
$$
\left(\frac{1}{1-t\frac{e^x-1}{x}}\right)^k=\sum\limits_{n\geqslant 0} \sum\limits_{m\geqslant 0} D(n,m,k)x^nt^m,
$$
by using the formula from \cite{Kru}, we have
$$
D(n,m,k)=\sum\limits_{i=0}^{n+m} T(n,m,i){i+k-1 \choose i}=\sum\limits_{k=0}^{n+m} \delta(m,i) \stirlingtwo{n+i}{i}\frac{i!}{(n+i)!}{i+k-1 \choose i}=
$$
$$
=\stirlingtwo{n+m}{m}\frac{m!}{(n+m)!}{m+k-1 \choose m}.
$$

According the Lagrange inverse theorem for power series $u(x,t)$ satisfied for the functional equation 
$$
u=x\,F(g,t)
$$
where $F(x,t)$ - power series with $F(0,0)\neq 0$
there hold
$$
[x^n]u(x,t)=\frac{k}{n}[x^{n-k}]F(x,t)^n.
$$
Applying that on our case $F(x,t)=\frac{1}{1-t\frac{e^x-1}{x}}$
we get
$$
u=\frac{x}{1-t\frac{e^{u}-1}{u}}.
$$
The solution of the equation is
$$
u(x,t)^k=\sum\limits_{n>0} \sum\limits_{m\geqslant 0} \frac{k}{n}D(n-k,m,n)x^n\,t^m
$$
where $D(n,m,k)$ is coefficients of $F(x,t)^k$.

Then
$$
D(n-k,m,n)=\stirlingtwo{n+m-k}{m}\frac{m!}{(n+m-k)!}{m+n-1 \choose m}
$$
For $k=1$ we get
$$
u(x,t)=\sum\limits_{n>0} \sum\limits_{m\geqslant 0} \frac{1}{n}\stirlingtwo{n+m-1}{m}\frac{m!}{(n+m-1)!}{m+n-1 \choose m}x^n\,t^m
$$
or after simplification we have
$$
u(x,t)=\sum\limits_{n>0} \sum\limits_{m\geqslant 0} \frac{1}{n!}\stirlingtwo{n+m-1}{m}x^n\,t^m.
$$

Hence, funtion  $u(x,t)$ is the compositional inverse function for $y(x,t)=(x-t(e^x-1))$.
The Lambert function is defined by \cite{LambertW}
$$
x=W(x)e^{W(x)}.
$$
It is well-known that for the equation
$$
g(x)=f(x)e^{f(x)}
$$
the solution is
$$
f(x)=W(g(x)).
$$
Next we apply that for our case 
$$
y=x(y,t)-t\,e^{x(t,y)}+t.
$$
Replacing
$$
z(y,t)=x(y,t)+t-y
$$
our equation will be
$$
z(y,t)=t\,e^{Z(y,t)-t+y}
$$
or
$$
z(y,t)e^{-Z(y,t)}=t\,e^{y-t}.
$$
Then the solution for $z(y,t)$ is
$$
z(y,t)=-W(-t\,e^{y-t})
$$
Hence,
$$
x(y,t)=y-t-W(-t\,e^{y-t})
$$
Considering the following product 
$$
E(x,t)=\sum_{n>0}P_n(t)\frac{x^n}{n!}\,(1-t)^{2n-1}=\frac{1}{1-t}\sum_{n>0}P_n(t)\frac{(x(1-t)^2)^n}{n!}
$$
we get
$$
E(x,t)=\frac{E_2(x(1-t)^2,t)}{(1-t)}=\frac{x(1-t)^2-t-W\left(-t\,e^{x(1-t)^2-t}\right)}{1-t}.
$$

differentiating with respect to $x$  the expression for $E(x,t)$ according properties of the derivative of the Lambert funtion $W(x)$ we arrive to the desired generating function
$$
\sum\limits_{n>0} \sum_{m>0} \eulereantwo{n}{m}\frac{x^n}{n!}t^m={{1-t}\over{W\left(-t\,e^{\left(1-t\right)^2\,x-t}
 \right)+1}}.
$$

\section{Identities}

First we prove an identity for the Stirling number of the second kind that presented in \cite{Knuth} (see the formula 6.4.3):

$$
u(x,t)=\sum\limits_{n>0} \sum\limits_{m\geqslant 0} \frac{1}{n!}\stirlingtwo{n+m-1}{m}\,t^m
$$
$$
E(x,t)=\frac{u(x(1-t)^2,t)}{1-t}=\sum\limits_{m\geqslant 0} \frac{1}{n!}\eulereantwo{n-1}{m}x^n\,t^m
$$
Then
$$
(1-t)E\left(\frac{x}{(1-t)^2},t\right)=\frac{1}{n!}\stirlingtwo{n+m-1}{m}\,t^m
$$
For the coefficients of composition $E(x,t)$ and $\frac{x}{(1-t)^2}$
we find coefficients for $\left(\frac{x}{(1-t)^2}\right)^k$
$$
\left(\frac{x}{(1-t)^2}\right)^k=\sum_{n>0}\sum_{m>0}p(n,m,k)=\delta(n,k){m+2k-1 \choose m}x^n\,t^m,
$$
where
$$
p(n,m,k)=\delta(n,k){m+2k-1 \choose m}.
$$
$$
E\left(\frac{x}{(1-t)^2},t\right)=\sum_{n}\sum_{m}e(n,m)x^n\,t^m,
$$
where
$$e(n,m)=\sum_{i=0}^{m}{\sum_{k=0}^{n+m-i}{\frac{1}{n!}\eulereantwo{k-1}{i}\,
 \delta_{k, n} \,{{m+2\,k-i-1}\choose{m-i}}}}=\sum_{i=0}^{m}{\frac{1}{n!}\eulereantwo{n-1}{i}\,
 {{m+2\,n-i-1}\choose{m-i}}}$$
Then coefficients for $(1-t)E\left(\frac{x}{(1-t)^2},t\right)$ are equal to
$$
e(n,m)-e(n,m-1)=\sum_{i=0}^{m}{\frac{1}{n!}\eulereantwo{n-1}{i}{{m+2\,n-i-1}\choose{m-i}}}-\sum_{i=0}^{m-1}{\frac{1}{n!}\eulereantwo{n-1}{i}\,
 {{m+2\,n-i-2}\choose{m-i-1}}}
$$
Note that inner expression of the second sum for $i=m$ is equal to 0. Then both sums could be combined
$$
e(n,m)-e(n,m-1)=\sum_{i=0}^{m}{\frac{1}{n!}\eulereantwo{n-1}{i}\left({{m+2\,n-i-1}\choose{m-i}}-{{m+2\,n-i-2}\choose{m-i-1}}\right)}=
$$
$$
=\sum_{i=0}^{m}{\frac{1}{n!}\eulereantwo{n-1}{i}{{m+2\,n-i-2}\choose{m-i}}}
$$
Then
$$
\frac{1}{n!}\stirlingtwo{n+m-1}{m}=\sum_{i=0}^{m}{\frac{1}{n!}\eulereantwo{n-1}{i}{{m+2\,n-i-2}\choose{m-i}}}
$$

$$
\stirlingtwo{n+m}{m}=\sum_{i=0}^{m}{\eulereantwo{n}{i}{{m+2\,n-i}\choose{m-i}}}
$$

Next we find new identities for the Euler numbers of the second kind based on the generating function $Eu(x,t)$. For that we consider the following composition of generating functions $Eu(x,t)$ and $\frac{t\,(x+1)}{(1-t)^2}$ 

$$
Eu\left(\frac{t\,(x+1)}{(1-t)^2},t\right)={{1-t}\over{W\left(-t\,e^{x\,t}
 \right)+1}} 
$$

$$\left(\frac{t\,(x+1)}{(1-t)^2}\right)^k=\sum_{n\geqslant 0}\sum_{m>0}{{k}\choose{n}}\,{{m+k-1}\choose{m-k}}x^nt^m
$$

\begin{equation}\label{Ident4}
Eu\left(\frac{t\,(x+1)}{(1-t)^2},t\right)=\sum_{n>0}\sum_{m>0}\sum_{i=0}^{m}{\sum_{k=0}^{n+m-i} \frac{1}{k!}\eulereantwo{k}{i}\,{{k}\choose{n}}\,{{m-i+k-1}\choose{m-i-k}}}x^n\,t^m
\end{equation}

For $x=0$ we have
$$
{{1-t}\over{W(-t)+1}} 
$$

For the Lambert function there hold
$$
\ln(W(x))=W(x)-\ln(x), x>0
$$
Differentiating, we get
$$
\frac{W'(x)}{W(x)}=\frac{1}{x}-W'(x)
$$
or
$$
W'(x)=\frac{W(x)}{x(1+W(x))}.
$$
Then
$$
\frac{1}{x}-W'(x)=\frac{1}{x(1+W(x))}
$$
$$
\frac{1}{(1+W(x))}=1-x\,W'(x)
$$
Substitute $x=-x$
$$
\frac{1}{(1+W(-x))}=1+x\,W'(-x)
$$
$$
W_0(-x)=\sum_{n>0} \frac{n^{n-1}}{n!}x^n
$$
$$
\frac{1}{(1+W(-x))}=\sum_{n\geqslant0} \frac{n^n}{n!}x^n
$$
$$
\frac{1-x}{(1+W(-x))}=1+\sum_{n>0} \left(\frac{n^n}{n!}-\frac{(n-1)^{n-1}}{(n-1)!}\right)x^n=1+\sum_{n>0} \left({n^{n-1}-(n-1)^{n-1}}\right)\frac{x^n}{(n-1)!}
$$

On the other side for $x=0$ the dual series will be the power series with one variable $t$ with $n=0$. The coefficients of the power series are defined by
$$
\sum_{i=0}^{m}{\sum_{k=0}^{m-i} \frac{1}{k!}\eulereantwo{k}{i}\,{{k}\choose{0}}\,{{m-i+k-1}\choose{m-i-k}}}
$$
Then we have the following identity
$$
\sum_{i=0}^{n}{\sum_{k=0}^{n-i} \frac{1}{k!}\eulereantwo{k}{i}\,{{n-i+k-1}\choose{n-i-k}}}=\frac{1}{(n-1)!}\left(n^{n-1}-(n-1)^{n-1}\right), n>0.
$$

For $x=-1$ we have the following identity

$$
{{1-t}\over{W(-t\,e^{-t})+1}}=1. 
$$

Considering (\ref{Ident4}), we get the explicit expression for $r(n,m)$
$$
r(n,m)=\sum_{i=0}^{m}{\sum_{k=0}^{n+m-i} \frac{1}{k!}\eulereantwo{k}{i}\,{{k}\choose{n}}\,{{m-i+k-1}\choose{m-i-k}}}.
$$

Note that $k\geqslant 0$, otherwise ${{k}\choose{n}}=0$. Then
$$
r(n,m)=\sum_{i=0}^{m}{\sum_{k=n}^{n+m-i} \frac{1}{k!}\eulereantwo{k}{i}\,{{k}\choose{n}}\,{{m-i+k-1}\choose{m-i-k}}}
$$

$$
r(n,m)=\sum_{i=0}^{m}{\sum_{k=0}^{m-i} \frac{1}{k!}\eulereantwo{k}{i}\,{{k+n}\choose{n}}\,{{m-i+k+n-1}\choose{m-i-k-n}}}
$$
Note that $r(n,m)$ is numeric triangle with $m\geqslant n$, because for $m<n$ there hold 
$$
{{m-i+k+n-1}\choose{m-i-k-n}}=0, 
$$

Consider (\ref{Ident4}) for $x=-1$:
$$
\sum_{m\geqslant 0}\sum_{n\geqslant 0} r(n,m)(-1)^n\,t^m
$$
Since $n<m$
$$
\sum_{m\geqslant 0}\sum_{n=0}^m r(n,m)(-1)^n\,t^m=1
$$
$$
\sum_{n=0}^m r(n,m)(-1)^n=0.
$$
we obtain
$$\sum_{m=0}^{n}{\left(-1\right)^{m}\,\sum_{i=0}^{n}{\sum_{k=0}^{i}{
 {{m+k}\choose{m}}\,\eulereantwo{m+k}{n-i}\,{{m+k+i-1
 }\choose{2\,m+2\,k-1}}}}}=0$$


\end{document}